\documentclass[final]{siamltex}
\usepackage{amsmath}
\usepackage{amssymb}
\usepackage{array}
\usepackage{blkarray}
\usepackage{capt-of}
\usepackage{color}
\usepackage{colortbl}
\usepackage{fancyhdr}
\usepackage{fancyvrb}
\usepackage{float}
\usepackage[T1]{fontenc}
\usepackage{framed}
\usepackage{geometry}
\usepackage{graphicx}
\usepackage{hyperref}
\usepackage[latin1]{inputenc}
\usepackage{stmaryrd}
\usepackage{subfigure}
\usepackage{upquote}
\usepackage{url}
\usepackage[usenames,dvipsnames,svgnames,table]{xcolor}
\usepackage{xspace}
\newcommand{\nc}{\newcommand}
\nc{\dsp}{\displaystyle}
\nc{\txt}{\textstyle}
\nc{\reff}[1]{(\ref{#1})}
\nc{\mrm}[1]{\mathrm{#1}}
\nc{\udl}[1]{\underline{#1}}
\nc{\ovl}[1]{\overline{#1}}
\nc{\al}{\underline{\boldsymbol{\alpha}}}
\nc{\la}{\underline{\boldsymbol{\lambda}}}
\nc{\llbr}{\llbracket}
\nc{\rrbr}{\rrbracket}
\nc{\lbr}{\lbrack}
\nc{\rbr}{\rbrack}
\nc{\N}{\mathbb{N}}
\nc{\Z}{\mathbb{Z}}
\nc{\D}{\mathbb{D}}
\nc{\Q}{\mathbb{Q}}
\nc{\R}{\mathbb{R}}
\nc{\C}{\mathbb{C}}
\nc{\T}{\mathbb{T}}
\nc{\StR}{\mathbb{S}^2_R}
\nc{\LtR}{\mathbb{L}^2_R}
\nc{\DtR}{\mathbb{D}^2_R}
\nc{\tld}[1]{\tilde{#1}}
\nc{\wtld}[1]{\widetilde{#1}}
\nc{\hu}{\hat{u}}
\nc{\wh}[1]{\widehat{#1}}
\nc{\sumeven}{\sum_{k=-N/2}^{N/2}{\hspace{-0.3cm}}'{\;\,}}
\nc{\sumodd}{\sum_{k=-\frac{N-1}{2}}^{\frac{N-1}{2}}}
\nc{\sumoddl}{\sum_{l=-\frac{N-1}{2}}^{\frac{N-1}{2}}}
\nc{\cqfd}{~\hbox{\vrule width 2.5pt depth 2.5 pt height 3.5 pt}}
\title{Computing hyperbolic choreographies}
\author{Hadrien Montanelli\thanks{Oxford University Mathematical Institute, Oxford OX2 6GG, UK.\ \ Supported by 
the European Research Council under the European Union's Seventh Framework Programme (FP7/2007--2013)/ERC grant agreement 
no.\ 291068.\ \ The views expressed in this article are not those of the ERC or the European Commission, and the European Union is not 
liable for any use that may be made of the information contained here.} }

\begin{document}

\maketitle

\begin{abstract}
An algorithm is presented for numerical computation of choreographies in spaces of constant negative curvature in a hyperbolic cotangent potential,
extending the ideas given in a companion paper~\cite{montanelli2016a} for computing choreographies
in the plane in a Newtonian potential and on a sphere in a cotangent potential.
Following an idea of Diacu, P\'{e}rez-Chavela and Reyes Victoria \cite{diacu2012d}, we use stereographic projection and study the problem in the Poincar\'{e}
disk.
Using approximation by trigonometric polynomials and optimization methods with exact gradient and exact Hessian matrix,
we find new choreographies, hyperbolic analogues of the ones presented in~\cite{montanelli2016a}.
The algorithm proceeds in two phases: first BFGS quasi-Newton iteration to get close to a solution, then Newton iteration for high accuracy.
\end{abstract}



\pagestyle{myheadings}
\thispagestyle{plain}

\markboth{COMPUTING HYPERBOLIC CHOREOGRAPHIES}{MONTANELLI}

\section{Introduction} 

Following the work of Chernoivan and Mamev~\cite{chernoivan1999} and Kilin~\cite{kilin1999}, there has been a growing interest in the $n$-body problem in spaces of constant curvature,
led by Borisov and his collaborators \cite{borisov2006,borisov2007,borisov2004}, Diacu and his collaborators~\cite{diacu2012c,diacu2013b,diacu2011,diacu2012d,diacu2012a,diacu2012b,perez2012}
and others \cite{carinena2005,shchepetilov2006}.
Recently, using numerical methods, the author has found new periodic solutions in the positive curvature case~\cite{montanelli2016a}
(i.e., on the sphere $\StR$ of radius $R$).
These are very special periodic configurations in which the bodies share a common orbit and are uniformly spread along it,
the \textit{spherical choreographies}.
Curved versions of the \textit{planar choreographies} found by Sim\'{o} in the early 2000s \cite{simo2001},
they can be computed to high accuracy using stereographic projection, trigonometric interpolation and optimization.
We show in this paper how these ideas can be used to find choreographies in spaces of negative curvature $-1/R^2$,
the \textit{hyperbolic choreographies}.
These are hyperbolic analogues of the planar and spherical choreographies and, as $R\rightarrow+\infty$, they converge to the planar
ones at a rate proportional to $1/R^2$.
 
\section{Hyperbolic choreographies}

While there really is only one model of two-dimensional spherical geometry (the sphere $\StR$ with the great-circle distance),
there are several models of hyperbolic geometry, 
including the Beltrami-Klein disk, the Poincar\'{e} disk, the Poincar\'{e} half-plane and the Lorentz hyperboloid models, with appropriate geodesic distances.
In this paper, we first use the latter and then reformulate the problem on the Poincar\'{e} disk using stereographic projection, following~\cite{diacu2012d}. 
Once on the disk, we use the techniques presented in~\cite{montanelli2016a}.

\begin{figure}
\hspace{-1.5cm}
\includegraphics [scale=.4]{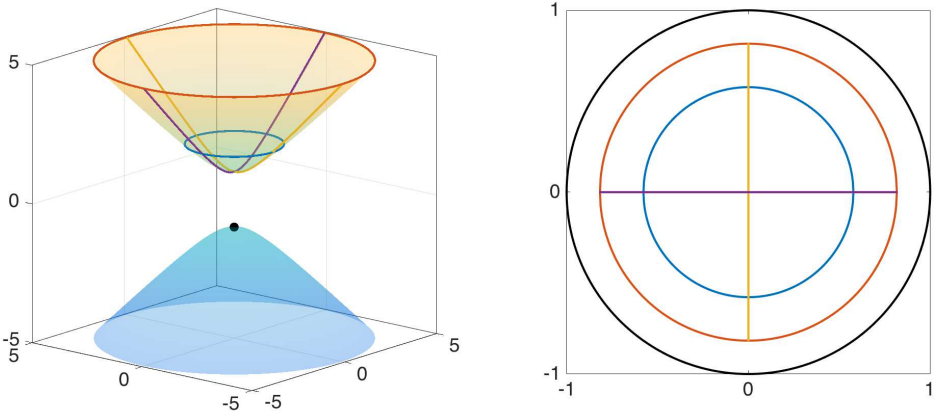}
\caption{\textit{On the left, a two-sheeted hyperboloid with $R=1$. The space $\LtR$ corresponds to the forward sheet only. 
The blue and the red curves are sections $x_3=2$ and $x_3=5$ and correspond to circles of equations $x_1^2+x_2^2=-R^2+2^2$ and  $x_1^2+x_2^2=-R^2+5^2$.
The yellow and the purple curves are sections $x_1=0$ and $x_2=0$ and correspond to hyperbolas of equations $x_2^2 - x_3^2=-R^2$ and $x_1^2 - x_3^2=-R^2$.
The north pole $(0,0,-R)$ of the backward sheet is indicated by a black dot.
On the right, the Poincar\'{e} disk $|z|<1$, obtained using stereographic projection from the north pole of the backward sheet. The curves correspond to the projections of 
the curves of the same colour on the forward sheet. The black line, corresponding to $|z|=1$, is the projection of the point at infinity.}}
\label{figure1}
\end{figure}

To describe hyperbolic geometry, the Lorentz model uses the forward sheet of a two-sheeted hyperboloid, defined as
\begin{equation}
\LtR = \{X=(x_1,x_2,x_3)^T\in \R^3, \; X\odot X=-R^2, \; x_3>0\}, \quad R>0,
\end{equation}

\noindent with Lorentz inner product
\begin{equation}
X\odot Y = x_1y_1 + x_2y_2 - x_3y_3
\end{equation}

\noindent and Lorentz distance
\begin{equation}
d(X,Y) = \sqrt{(X-Y)\odot(X-Y)}
\label{Lorentzmetric}
\end{equation}

\noindent for $X=(x_1,x_2,x_3)^T$ and $Y=(y_1,y_2,y_3)^T$ on $\LtR$. 
A two-sheeted hyperboloid with $R=1$ is shown in Figure~\ref{figure1}.
The geodesic distance between $X$ and $Y$ on $\LtR$ is
\begin{equation}
\hat{d}(X,Y) = R\,\mrm{acosh}\Big(-\frac{X\odot Y}{R^2}\Big).
\label{geodesicdistance}
\end{equation}

\noindent The forward sheet $\LtR$ together with the geodesic distance \reff{geodesicdistance} is called
the \textit{Lorentz hyperboloid model}.
Note that this model uses \textit{extrinsic coordinates}: $\LtR$ is embedded in $\R^3$, i.e., points on $\LtR$
are represented by Cartesian coordinates in $\R^3$.

The $n$-body problem on $\LtR$ describes the motion of $n$ bodies on $\LtR$ with Cartesian coordinates 
$X_j(t)\in\R^3$, $0\leq j\leq n-1$, via the $n$ coupled nonlinear ODEs
\begin{equation}
\dsp X_j''(t) - \sum_{\underset{i\neq j}{i=0}}^{n-1}\frac{R^3X_i(t) + R(X_i(t)\odot X_j(t))X_j(t)}{\big[(X_i(t)\odot X_j(t))^2-R^4\big]^{3/2}}
- R^{-2}\big(X'_j(t)\odot X'_j(t)\big) X_j(t) = 0, \; 0\leq j\leq n-1.
\label{hypnewton}
\end{equation}

\noindent The potential associated with \reff{hypnewton} is a hyperbolic cotangent potential. 
It is a generalization of the Newtonian potential and dates back to the 19th century with the work of Bolyai, Lobachevsky and Killing.

We are looking for \textit{hyperbolic choreographies}, i.e., solutions $X_j(t)$ such that
\begin{equation}
X_j(t) = Q\Big(t + \frac{2\pi j}{n}	\Big), \quad 0\leq j\leq n-1,
\label{hypchoreographies}
\end{equation}

\noindent for some $2\pi$-periodic function $Q:[0,2\pi]\rightarrow\LtR$. 
We can choose the period equal to $2\pi$ since if $Q(t)$ is a $T$-periodic of \reff{hypnewton} on $\LtR$
then $\lambda^{-2/3}Q(\lambda t)$, $\lambda=T/(2\pi)$, is a $2\pi$-periodic solution in $\mathbb{L}^2_{R'}$ 
with $R'=\lambda^{-2/3}R$.
As in the plane and on the sphere, they correspond to minima of the action associated with \reff{hypnewton}, defined
as the integral over one period of the kinetic minus potential energy, 
\begin{equation}
A = \int_0^{2\pi}\big(K(t) - U(t)\big)dt,
\end{equation}

\noindent with kinetic energy
\begin{equation}
\dsp K(t) = \frac{1}{2}\sum_{j=0}^{n-1} X'_j(t)\odot X'_j(t) = \frac{1}{2}\sum_{j=0}^{n-1} 
Q'\Big(t + \frac{2\pi j}{n}	\Big)\odot Q'\Big(t + \frac{2\pi j}{n}	\Big)
\label{Kinetic}
\end{equation}

\noindent and potential energy
\begin{equation}
\dsp U(t) = -\frac{1}{R}\sum_{j=0}^{n-1}\sum_{i=0}^{j-1} \coth \frac{\hat{d}(X_i(t),X_j(t))}{R}.
\label{hypcotangentpotential}
\end{equation}

\noindent Using the trigonometric identity $\cot(\mrm{acosh}(x))=x/\sqrt{x^2-1}$, the potential energy can be rewritten
\begin{equation}
\dsp U(t) = \frac{1}{R}\sum_{j=0}^{n-1}\sum_{i=0}^{j-1}\frac{X_i(t)\cdot X_j(t)}{\sqrt{(X_i(t)\odot X_j(t))^2 - R^4}}.
\label{Potential}
\end{equation}

\noindent Since the integral of \reff{Kinetic} does not depend on $j$ and the integral of \reff{Potential} only depends
on $i-j$, the action is given by
\begin{equation}
\dsp A = \frac{n}{2}\int_0^{2\pi} Q'(t)\odot Q'(t) \, dt
- \frac{n}{2R}\sum_{j=1}^{n-1}\int_0^{2\pi} \frac{Q(t)\odot Q\big(t + \frac{2\pi j}{n}\big)}{\sqrt{\big(Q(t)\odot Q\big(t + \frac{2\pi j}{n}\big)\big)^2 - R^4}}\,dt.
\label{hypaction}
\end{equation}

\noindent We are also looking for \textit{relative} hyperbolic choreographies, 
\begin{equation}
X_j(t) = R_\omega(t) \, Q\Big(t + \frac{2\pi j}{n}	\Big), \quad 0\leq j\leq n-1, 
\quad R_\omega(t)=
\begin{bmatrix}
\cos(\omega t) & -\sin(\omega t) & 0 \\
\sin(\omega t) & \cos(\omega t) & 0 \\
0 & 0 & 1
\end{bmatrix},
\label{hypchoreographies}
\end{equation}

\noindent i.e., choreographies rotating with angular velocity $\omega$ along the $x_3$-axis. 
In this case, the kinetic part of \reff{hypaction} is
\begin{equation}
\frac{n}{2}\int_0^{2\pi} \big(R_\omega(t) Q'(t) + R'_\omega(t) Q(t)\big)\odot \big(R_\omega(t) Q'(t) + R'_\omega(t) Q(t)\big) \, dt.
\label{hypaction2}
\end{equation}

\noindent Hyperbolic choreographies correspond to functions $Q(t)$ which minimize \reff{hypaction}--\reff{hypaction2}.

Now, let us reformulate this minimization problem on the Poincar\'{e} disk using stereographic projection.
Points $X=(x_1,x_2,x_3)^T$ on $\LtR$ are mapped to points $z=P_R(X)$ on the Poincar\'{e} disk $\DtR=\{z\in\C, \, \vert z\vert<R\}$
via
\begin{equation}
z = P_R(X) = \frac{R x_1 + i R x_2}{R + x_3}.
\label{stereoproj}
\end{equation}

\noindent The inverse mapping is given by
\begin{equation}
X = P_R^{-1}(z) = \frac{1}{R^2-\vert z\vert^2}(2R^2\mrm{Re}(z),2R^2\mrm{Im}(z),R^3 + R\vert z\vert^2)^T.
\label{stereoprojinv}
\end{equation}

\noindent Note that \reff{stereoproj} is a stereographic projection from the north pole $(0,0,-R)$ of the 
backward sheet of the hyperboloid---see Figure~\ref{figure1} for an example of such a projection.
The Lorentz distance \reff{Lorentzmetric} between two points on $\LtR$ is transformed into the distance $d(z,\xi)$
between their projections $z=P_R(X)$ and $\xi=P_R(Y)$ defined as
\begin{equation}
d(z,\xi) = \frac{2R^2\vert z - \xi \vert}{\sqrt{(R^2-\vert z\vert^2)(R^2-\vert \xi\vert^2)}},
\label{Lorentzmetric2}
\end{equation}

\noindent and the geodesic distance \reff{geodesicdistance} into
\begin{equation}
\hat{d}(z,\xi) = 2R\,\mrm{asinh}\frac{d(z,\xi)}{2R}.
\label{geodesicdistance2}
\end{equation}

\noindent The Poincar\'{e} disk $\DtR$ together with the geodesic distance \reff{geodesicdistance2} 
is called the \textit{Poincar\'{e} disk model}.
This model uses \textit{intrinsic coordinates} since points on $\DtR$ are represented by complex coordinates, i.e.,
$\DtR$ is not embedded in any higher dimensional space.
Let $q(t)=P_R(Q(t))$ denote the projection of $Q(t)\in\LtR$ onto $\DtR$, and
\begin{equation}
z_j(t) = P_R(X_j(t)) = P_R\Big(Q\Big(t + \frac{2\pi j}{n}	\Big)\Big) = q\Big(t + \frac{2\pi j}{n}	\Big), \quad 0\leq j\leq n-1,
\end{equation}

\noindent the projections of the $n$ bodies. The kinetic part \reff{hypaction2} of the action can be rewritten as 
\begin{equation}
\dsp\frac{n}{2}\int_0^{2\pi}\bigg(\frac{2R^2\vert q'(t) + i\omega q(t)\vert}{R^2 - \vert q(t)\vert^2}\bigg)^2dt
\end{equation}

\noindent with conformal factor $4R^4/(R^2 - \vert q(t)\vert^2)^2$. 
To derive the formula for the potential part of \reff{hypaction} in intrinsic coordinates, let us come back to the potential energy \reff{Potential}.
On the Poincar\'{e} disk $\DtR$, it is given by
\begin{equation}
\dsp U(t) = -\frac{1}{R}\sum_{j=0}^{n-1}\sum_{i=0}^{j-1} \coth \frac{\hat{d}(z_i(t),z_j(t))}{R} = 
-\frac{1}{R}\sum_{j=0}^{n-1}\sum_{i=0}^{j-1} \coth \bigg(2\,\mrm{asinh}\frac{d(z,\xi)}{2R}\bigg).
\end{equation}

\noindent Using the trigonometric identity $\coth(2\,\mrm{asinh}(x/2))=(x^2/2+1)/(x\sqrt{x^2/4+1})$ and integrating over one period,
we find that the action is given by
\begin{equation}
\dsp A =  \dsp\frac{n}{2}\int_0^{2\pi}\bigg(\frac{2R^2\vert q'(t) + i\omega q(t)\vert}{R^2 - \vert q(t)\vert^2}\bigg)^2dt
+ \frac{n}{2R}\sum_{j=1}^{n-1}\int_0^{2\pi}\frac{2R^2 + D_j(t)^2}{D_j(t)\sqrt{4R^2 + D_j(t)^2}}\,dt,
\label{hypaction3}
\end{equation}

\noindent with $D_j(t) = d\big(q(t),q\big(t+\frac{2\pi j}{n}\big)\big)$. 
Hyperbolic choreographies correspond to functions $q(t)$ which minimize \reff{hypaction3}.

\section{Computing hyperbolic choreographies}

Our method for computing hyperbolic choreographies is based on the algorithm presented in \cite{montanelli2016a}---we summarize here quickly the key ideas behind this algorithm, 
and refer to \cite{montanelli2016a} for details. 

The algorithm uses trigonometric interpolation and numerical optimization of the action \reff{hypaction3}. The function $q(t)$ is represented by its trigonometric interpolant in the $\exp(ikt)$ basis.
The optimization variables are the real and imaginary parts of its Fourier coefficients and the action is computed with the exponentially accurate trapezoidal rule \cite{trefethen2014}. 
Closed-form expressions for the gradient and the Hessian of the action with respect to the optimization variables can be derived and are used in the numerical optimization, 
which is carried out in two phases.

{\em Phase 1. Quasi-Newton optimization methods.} Numerical optimization methods with the exact gradient and based on approximations of the Hessian are employed with a small number of optimization variables.
The accuracy of the solution at this stage is from one to five digits.
This phase is computationally very cheap. 

{\em Phase 2. Newton's method.} Once an approximation to a choreography has been computed via a quasi-Newton method, one can improve the accuracy to typically ten digits with a few steps of Newton's method with exact Hessian, and a larger number of optimization variables.
This phase is computationally more expensive.

We use Chebfun \cite{chebfun}, its extension to periodic problems \cite{montanelli2015b} and MATLAB \texttt{fminunc} code for our computations. 
Once a choreography has been computed by our algorithm, we check that its Fourier coefficients decay to sufficiently small values, the gradient of the action \reff{hypaction3} has small norm and 
that it is a solution of the equations of motion \reff{hypnewton} projected onto the Poincar\'{e} disk. The latter were first derived in \cite{diacu2012d} and are given by
\begin{equation}
z_j''(t) = -\frac{2\bar{z}_j(t)z_j'^2(t)}{R^2-\vert z_j(t)\vert^2} + \frac{4R}{\lambda_j(t)}\sum_{\underset{i\neq j}{i=0}}^{n-1}\frac{P_{j,i}(t)}{\Theta_{j,i}(t)^{3/2}}, \quad 0\leq j\leq n-1,
\label{hypnewtonbis}
\end{equation}

\noindent where $\lambda_j(t)=4R^4/(R^2-\vert z_j(t)\vert^2)^2$ is the conformal factor introduced before, while $P_{j,i}(t)$ and $\Theta_{j,i}(t)$ are defined by
\begin{equation}
P_{j,i}(t) = \big[R^2-\vert z_j(t)\vert^2\big]\big[R^2-\vert z_i(t)\vert^2\big]^2\big[R^2-\bar{z}_i(t)z_j(t)\big]\big[z_i(t)-z_j(t)\big],
\end{equation}

\noindent and 
\begin{equation}
\begin{array}{ll}
\Theta_{j,i}(t) & = \big[2R^2z_j(t)\bar{z}_i(t) + 2R^2z_i(t)\bar{z}_j(t) - (\vert z_j(t)\vert^2+R^2)(\vert z_i(t)\vert^2+R^2)\big]^2\\
& - \big[R^2-\vert z_j(t)\vert^2\big]^2\big[R^2-\vert z_i(t)\vert^2\big]^2.
\end{array}
\end{equation}

\begin{figure}
\hspace{-1.6cm}
\includegraphics [scale=.53]{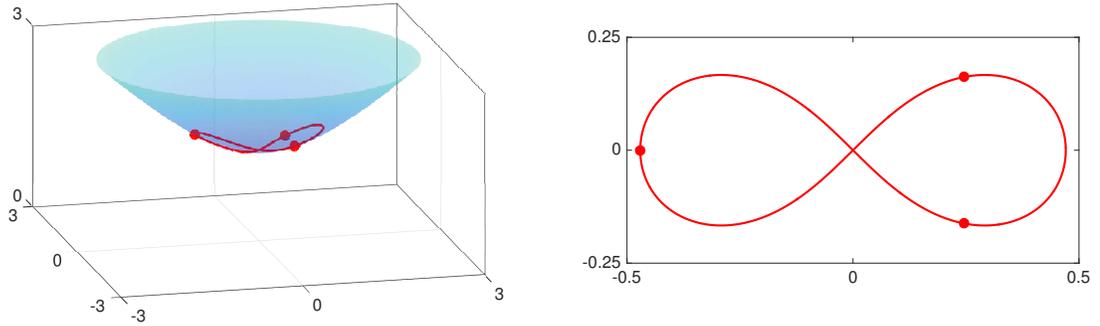}
\caption{\textit{Hyperbolic figure-eight with $R=1.5$ (left) and its projection on the Poincar\'{e} disk (right). The dots show the bodies at time $t=0$.
This choreography can be computed to about twelve digits of accuracy in less than $1.5$ seconds on a $2.7$\,{\normalfont GHz} Intel i$7$ machine.}}
\label{figure2}
\end{figure}

\begin{figure}
\centering
\includegraphics [scale=.35]{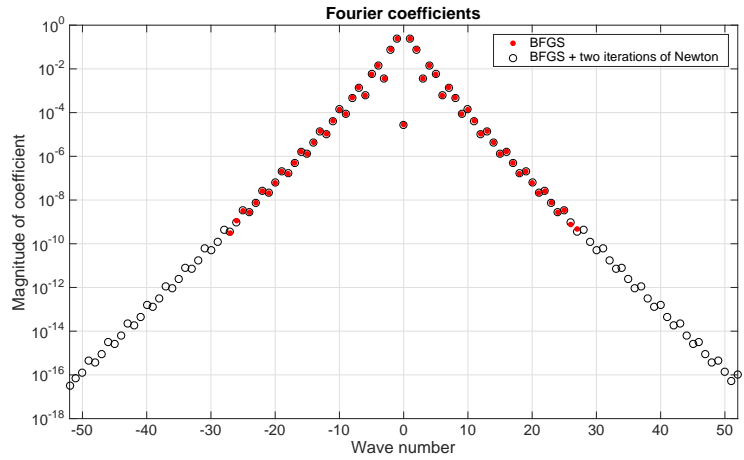}
\caption{\textit{Absolute values of the Fourier coefficients of the hyperbolic figure-eight of Figure~$\ref{figure2}$, obtained by BFGS (red dots) and BFGS followed by two steps of Newton's method (black circles).}}
\label{figure2bis}
\end{figure}

\begin{table}
\vspace{.4cm}
\begin{center}
\begin{tabular}{l|cc}
& Phase 1: BFGS & Phase 2: Newton\\
\hline
Action & 27.840867421590943 & 27.840867421590929\\
Number of coefficients & 55 & 105\\
Computer time (s) & 0.7734 & 0.6427\\
Number of iterations & 87 & 2\\
Relative $2$-norm of the gradient & 8.08e-08 & 3.87e-15\\
Smallest coefficient & 4.75e-10 & 1.03e-16\\
Relative $2$-norm of the residual & 9.38e-07 & 3.54e-13\\
\end{tabular}
\end{center}
\vspace{.2cm}
\caption{\textit{Two-phase computation of the hyperbolic figure-eight of Figure~$\ref{figure2}$.}}
\label{table:figure2}
\end{table}

The first choreography that we present is the hyperbolic figure-eight of the three-body problem with $R=1.5$, see Figure~\ref{figure2}.
Table~\ref{table:figure2} shows that, after 87 iterations of the first phase, the choreography satisfies \reff{hypnewtonbis} to six digits and after two iterations of the second phase it satisfies it to twelve digits.
Figure~\ref{figure2bis} shows the Fourier coefficients of the solution, they decay to about $10^{-10}$ after the first phase and to about $10^{-16}$ after the second phase.
We see in Table~\ref{table:figure2} that this choreography satisfies \reff{hypnewtonbis} to $12$ digits of accuracy after the second phase.

\begin{figure}
\hspace{-1.7cm}
\includegraphics [scale=.58]{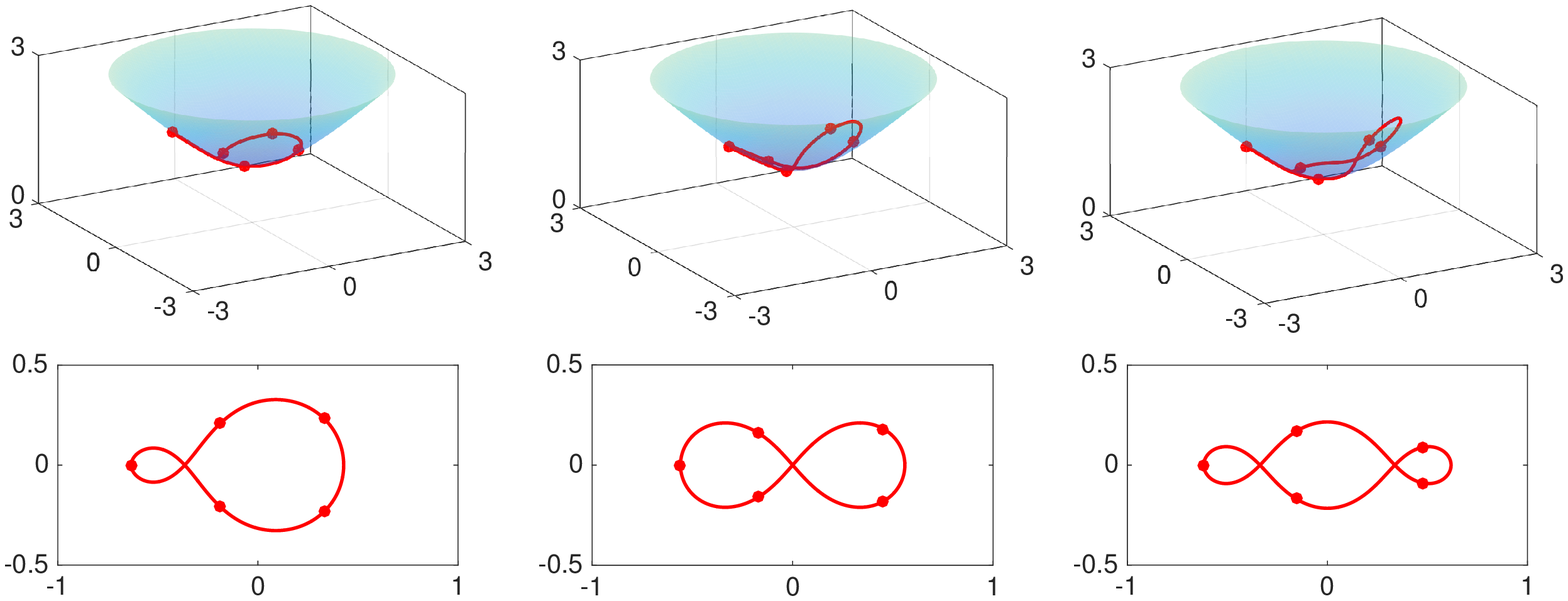}
\caption{\textit{Hyperbolic choreographies with $R=1.2$ (top) and their projections on the Poincar\'{e} disk (bottom). The dots show the bodies at time $t=0$.}}
\label{figure3}
\end{figure}

Many choreographies can be found with our algorithm. We show in Figure~\ref{figure3} three hyperbolic choreographies of the five-body problem with $R=1.2$. These are curved
versions of the choreographies found by Sim\'{o} in~\cite{simo2001}. As shown in Table~\ref{table:figure3}, they can be computed to high accuracy with a few hundred Fourier coefficients.

\begin{figure}
\hspace{-1.9cm}
\includegraphics [scale=.58]{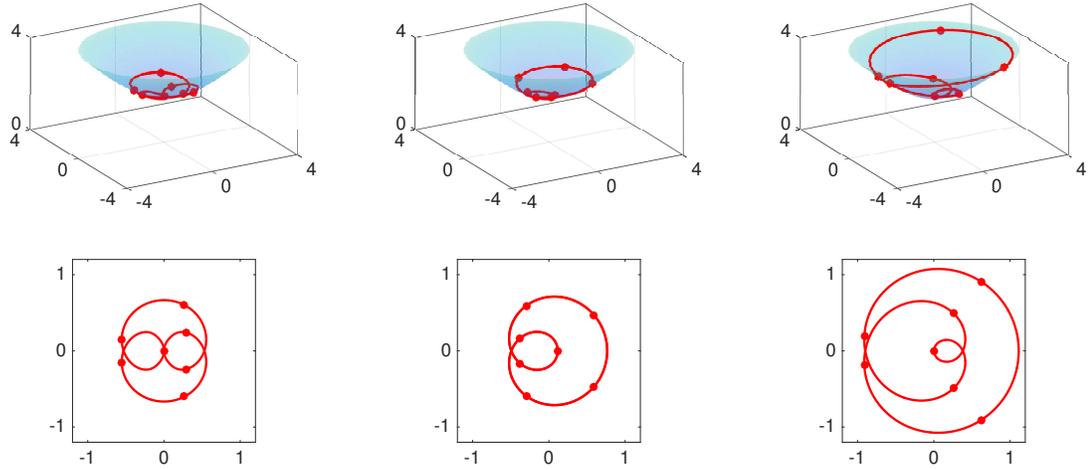}
\caption{\textit{First row: relative hyperbolic choreographies with $R=2$ and angular velocities $2.8$ (left), $-2.9$~(center) and $2.31$ (right). 
Second row: their projections on the Poincar\'{e} disk. The dots show the bodies at time $t=0$.}}
\label{figure4}
\end{figure}

Relative choreographies can also be computed, see Figure~\ref{figure4}. Again, a few hundred coefficients is enough to get about $10$-digit accuracy.

\begin{table}
\begin{center}
\begin{tabular}{l|cc|cc|cc}
& Phase 1 & Phase 2 & Phase 1 & Phase 2 & Phase 1 & Phase 2\\
\hline
Action & 88.8733 & 88.8733 & 90.6073 & 90.6073 & 96.2604 &  96.2604\\
Number of coefficients & 75 &  305 & 55 & 155 & 65 & 245\\
Computer time (s) & 2.52 & 17.54 & 0.75 & 3.02 & 1.03 & 10.88\\
Number of iterations & 112 & 4 & 68 & 2 & 105 & 4\\
Relative $2$-norm of the gradient & 2.76e-08 & 5.82e-13 & 1.03e-07 & 8.30e-15 & 2.96e-09 & 9.89e-15\\
Smallest coefficient                      & 1.55e-06 & 1.33e-17 & 1.13e-09 & 3.45e-18 & 3.85e-08 & 5.91e-18\\
Relative $2$-norm of the residual  & 3.56e-03 & 3.54e-12 & 7.23e-06 & 3.57e-13 & 4.58e-03 & 2.03e-12\\
\end{tabular}
\end{center}
\vspace{.2cm}
\caption{\textit{Computation of the hyperbolic choreographies of Figure~$\ref{figure3}$.}}
\label{table:figure3}
\end{table}

\section{Limit of infinitely large $R$}

\begin{figure}
\hspace{-.3cm}
\includegraphics [scale=.4]{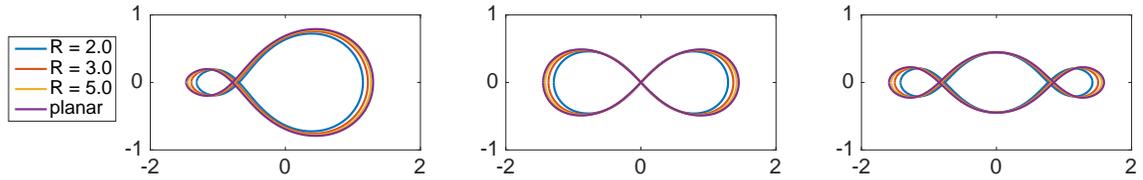}
\caption{\textit{Hyperbolic choreographies of Figure~$\ref{figure3}$ (multiplied by a factor of $2$) for different values of $R$. 
As $R$ increases, the hyperbolic choreographies converge to the planar ones.}}
\label{figure5}
\end{figure}

\begin{figure}
\centering
\includegraphics [scale=.35]{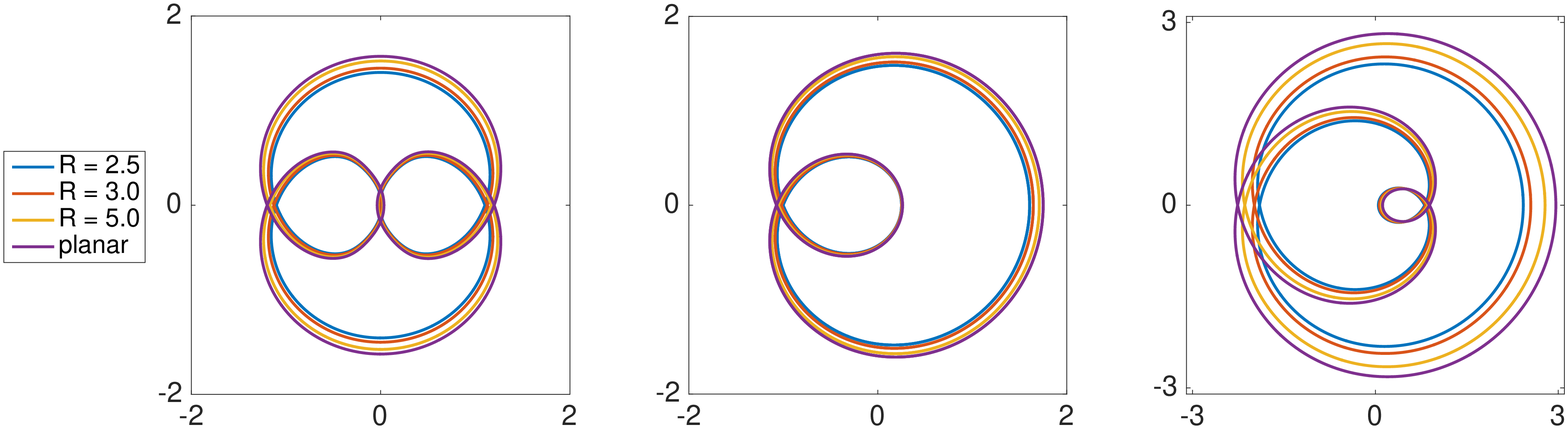}
\caption{\textit{Relative hyperbolic choreographies of Figure~$\ref{figure4}$ (multiplied by a factor of $2$) for different values of $R$.
As $R$ increases, the hyperbolic choreographies converge to the planar ones.}}
\label{figure6}
\end{figure}

In the limit $R\rightarrow\infty$, the Poincar\'e disk converges to the complex plane. The distances \reff{Lorentzmetric2} and \reff{geodesicdistance2}
converge to twice the absolute value and the action \reff{hypaction2} converges to four times the action in the plane, since it involves squares of distances.
As a consequence, twice the hyperbolic choreographies converge to the planar ones as $R\rightarrow\infty$, as shown in Figures \ref{figure5} and \ref{figure6}.
Tables \ref{table:figure5} and \ref{table:figure6} report the $\infty$-norm of the difference between analogous hyperbolic and planar choreographies as 
$R$ increases. The convergence appears to be at rate proportional to the absolute value of the curvature $1/R^2$.

\begin{table}
\begin{center}
\begin{tabular}{l|cccccc}
& $R$ = 2 & 3 & 5 & 10 & 100 & 1000\\
\hline
Left & 1.56e-01 & 7.74e-02 & 3.03e-02 & 7.87e-03 & 7.98e-05 &  7.99e-07\\
Middle & 1.59e-01 & 7.97e-02 & 3.09e-02 & 7.98e-03 &  8.07e-05 & 8.14e-07\\
Right & 1.70e-01 & 8.54e-02 & 3.31e-02 & 8.55e-03 & 8.65e-05 & 8.65e-07\\
\end{tabular}
\vspace{.2cm}
\caption{\textit{Convergence of the hyperbolic choreographies to the planar ones in Figure~$\ref{figure5}$.}}
\label{table:figure5}
\end{center}
\end{table}

\begin{table}
\begin{center}
\begin{tabular}{l|cccccc}
& $R$ = 2.5 & 3 & 5 & 10 & 100 & 1000\\
\hline
Left & 1.70e-01 & 1.25e-01 & 4.93e-02 & 1.28e-02 & 1.30e-04 &  1.31e-06\\
Middle & 1.43e-01 & 1.05e-01 & 4.07e-02 & 1.05e-02 &  1.07e-04 & 1.09e-06\\
Right & 5.34e-01 & 4.11e-01 & 1.77e-01 & 4.87e-02 & 5.03e-04 & 5.04e-06\\
\end{tabular}
\vspace{.2cm}
\caption{\textit{Convergence of the relative hyperbolic choreographies to the planar ones in Figure~$\ref{figure6}$.}}
\label{table:figure6}
\end{center}
\end{table}

\section{Discussion}

We have shown numerical evidence that choreographies also exist in spaces of constant negative curvature. As in the plane and on the sphere, they can be computed to high
accuracy using trigonometric interpolation and minimization of the action. 

The author believes that the techniques described in this paper can be applied not only to particle dynamics but also to other types of dynamics.
A possible extension of this work would therefore be the study of choreographies of the $n$-vortex problem~\cite{newton2001}, which
describes the motion of $n$ vortices, complex potentials associated with the two-dimensional, irrotational and incompressible Euler equations.

\bibliographystyle{siam}
\bibliography{Montanelli2016c}

\end{document}